\nonstopmode

\documentclass[11pt]{article}

\usepackage[fleqn]{amsmath}

\usepackage{hyperref}
\hypersetup{
    colorlinks=true,
    linkcolor=blue,
    filecolor=magenta,      
    urlcolor=blue,
}

\usepackage{tcolorbox}   
\definecolor{mycolor}{rgb}{0.122, 0.435, 0.698}

  \usepackage{geometry}
\usepackage{amsmath,amssymb,amsthm}
\usepackage{graphics,epsfig,calc}

\usepackage{latexsym,epsfig,bm,amssymb}
\usepackage{xcolor}
\usepackage{amsthm,mathrsfs}

\DeclareSymbolFont{AMSb}{U}{msb}m{n}
\DeclareSymbolFontAlphabet{\mathbb}{AMSb}

\newcommand{\beqn}{\begin{eqnarray}}
\newcommand{\eeqn}{\end{eqnarray}}
\newcommand{\be}{\begin{equation}}
\newcommand{\ee}{\end{equation}}

\newcommand{\bsp}{\begin{split}}
\newcommand{\esp}{\end{split}}

\newcommand{\ba}{\begin{array}}
\newcommand{\ea}{\end{array}}

\newcommand{\bpr}{\begin{proof}}
\newcommand{\epr}{\end{proof}}

\newcommand{\HS}{{\rm HS}}

\newcommand{\aA}{{\mathbb A}}

\newcommand{\cD}{{\cal D}}

\newcommand{\cT}{{\cal T}}

\newcommand{\da}{\dagger}
\newcommand{\ci}{\cite}
\newcommand{\de}{\delta}
\newcommand{{\De}}{{\Delta}}

\newcommand{\fr}{\frac}

\newcommand{\ga}{\gamma}

\newcommand{\la}{\label}

\newcommand{  \om}{  \omega}

\newcommand{ \ov}{ \overline}

\newcommand{\si}{\sigma}
\newcommand{\ti}{\tilde}

\newcommand{\ve}{\varepsilon}

\newcommand\C{{\mathbb C}}
\newcommand\R{{\mathbb R}}

\newcommand{\1}{{\hspace{0.1mm}}}

\newcommand{\tr}{\mathop{\rm tr\,}\nolimits}

 \newcommand{\st}{\stackrel}

   \newcommand{\toC}{\st{C(0,\infty)}{-\!\!\!-\!\!\!-\!\!\!\!\longrightarrow}}

\newtheorem{theorem}{Theorem}[section]

\newtheorem{definition}[theorem]{Definition}

\newtheorem{lemma}[theorem]{Lemma}
\newtheorem{example}[theorem]{Example}

\newtheorem{remark}[theorem]{Remark}

\newtheorem{remarks}[theorem]{Remarks}
\newtheorem{cor}[theorem]{Corollary}
\newtheorem{proposition}[theorem]{Proposition}

\newcommand{\bd}{\begin{definition}}
 \newcommand{\ed}{\end{definition}}

\newcommand{\bt}{\begin{theorem}}
 \newcommand{\et}{\end{theorem}}
\newcommand{\bqt}{\begin{qtheorem}}
 \newcommand{\eqt}{\end{qtheorem}}

\newcommand{\bp}{\begin{proposition}}
 \newcommand{\ep}{\end{proposition}}

\newcommand{\bl}{\begin{lemma}}
 \newcommand{\el}{\end{lemma}}
\newcommand{\bc}{\begin{cor}}
 \newcommand{\ec}{\end{cor}}

\newcommand{\bex}{\begin{example}}
 \newcommand{\eex}{\end{example}}
 
\newcommand{\bexs}{\begin{examples}}
 \newcommand{\eexs}{\end{examples}}

\newcommand{\bexe}{\begin{exercice}}
 \newcommand{\eexe}{\end{exercice}}

\newcommand{\br}{\begin{remark}}
 \newcommand{\er}{\end{remark}}
 
\newcommand{\brs}{\begin{remarks}}
 \newcommand{\ers}{\end{remarks}}

\newcommand{\bce}{\begin{center}}
\newcommand{\ece}{\end{center}}

\begin{document}
\begin{center}

{\huge On global dynamics for 
damped driven 

Jaynes--Cummings equations

}

\bigskip\smallskip

 {\large A.I. Komech$^1$ and E.A. Kopylova}\footnote{ 
 The research was funded in whole by Austrian Science Fund (FWF) 10.55776/PAT3476224.}
 \\
{\it
Institute of Mathematics of
BOKU
University, Vienna, Austria\\
}
 alexander.komech@boku.ac.at,\qquad
 elena.kopylova@boku.ac.at

\smallskip

\end{center}

\setcounter{page}{1}
\thispagestyle{empty}

\begin{abstract}

 The article concerns
  damped driven Jaynes--Cummings equation which
describes quantised one-mode Maxwell field
coupled to a  two-level 
molecule. We consider a broad class of  
 damping and 
 pumping which are polynomial in the creation and  annihilation operators,
 and  their structures
 correspond to
 the theory of
completely positive and trace preserving generators (CPTP) of Lindblad and Kossakowski \& al.

 Our
  main result is the construction of global generalised solutions with values   in the Hilbert
space of nonnegative  Hermitian Hilbert--Schmidt operators
  in the case of 
  time-dependent pumping.
  The proofs rely on finite-dimensional approximations of the annihilation and creation ope\-rators. 
  
         \end{abstract}

  \smallskip
  

\tableofcontents

\setcounter{equation}{0}
\section{Introduction}
The Jaynes--Cummings equation is the basic model of Quantum Optics,
and it is used for description of various aspects of the laser action. 
The model without damping and pumping was introduced in  
\ci{JC1963} (the survey can be found in \ci{LM2021}). 
Various versions of the pumping are considered in \ci{AGC1992,BHPSM2024,DKM1994,JA1993,SF2002}.
The damping was introduced analysing quantum spontaneous emission \ci{A1973}--\ci{A1974}, \ci{BR1997,BJ2007,SF2002,VW2006}.
We construct global solutions for all initial values from the space of Hilbert--Schmidt 
operators in the case of time-dependent pumping.
\smallskip

Denote 
$X=F\otimes\C^2$, where $F$ is the single-particle Hilbert space
endowed with 
an orthonormal basis $|n\rangle$, $n=0,1,\dots$,
and the corresponding annihilation and   creation operators $a$ and $a^\dag$:
\be\la{aa} 
 a|n\rangle=\sqrt{n}|n-1\rangle,\qquad 
 a^\dag|n\rangle=\sqrt{n+1}|n+1\rangle, \qquad [a,a^\da]=1.
\ee
We will consider  a  damped driven version of  the Jaynes--Cummings equation 
\be\la{JC}
\dot \rho(t)=\aA\rho(t):=-i[H(t),\rho(t)]+\ga D(t)\rho(t),\qquad t\ge 0,
\ee
where
the density operator  $\rho(t)$ of the coupled field-molecule system 
 is  a  Hermitian operator in $X$.
 The Hamiltonian  $H(t)$ is the sum 
\be\la{JCh}
H(t)=H_0+pH_1(t),\qquad H_0:=\om_c a^\da a+\fr12  \om_a\si_3,\quad H_1(t)= (a+a^\dag)\otimes \si_1+A^e(t).
\ee
Here $H_0$ is the Hamiltonian of the free field and atom without interaction,
while $pH_1(t)$ is the interaction Hamiltonian,
$\om_c>0$ is the cavity resonance frequency,  
$\om_a>0$ is the atomic 
frequency, and
  $p\in\R$ is proportional to the molecular dipole moment.
   The pumping is represented by  selfadjoint operators
   $A^e(t)$, 
   and
$\si_1$ and $\si_3$ are the Pauli matrices
 acting on the factor $\C^2$ in $F\otimes\C^2 $, so $[a,\si_k]=[a^\dag,\si_k]=0$.
Finally, $\ga>0$, and $D(t)$ is a dissipation operator. 
We will consider the operator
\be\la{tiGa}
D_1\rho=a\rho a^\dag-\fr12 a^\dag a\rho-\fr12 \rho a^\dag a, 
 \ee
 used in \cite{A1973}--\cite{A1974}, \cite{BR1997,  BJ2007, L1976, SF2002,VW2006}, 
 and also its suitable modifications.
%
  \bd
  $\HS$ is the Hilbert space of Hermitian 
    Hilbert--Schmidt operators with the inner product {\rm \ci{RS1980}} 
\be\la{HS}
\langle\rho_1,\rho_2\rangle_\HS=\tr[\rho_1\rho_2].
\ee
\ed
  \bd
i)  $|n,s_\pm\rangle=|n\rangle\otimes s_\pm$ is the orthonormal basis in $X$,
$s_\pm\in \C^2$ and $\si_3 s_\pm=\pm s_\pm$.
\smallskip\\
ii) $X_\infty$ is the space 
of all finite linear combinations of the vectors $|n,s_\pm\rangle$.
\smallskip\\
iii) $\cD\subset\HS$ is the subspace of finite rank Hermitian operators 
\be\la{fr1}
\rho=\sum_{n,n'=0}^\infty \,\sum_{s,s'=s_\pm}\rho_{n,s;n',s'}
|n,s\rangle\otimes \langle n',s'|.
\ee
\ed

      Our main goal is to construct global nonnegative solutions for the equation (\ref{JC})
    in the  Hilbert space $\HS$  in the case of time-dependent
pumping $A^e$.  The main issue is that
the operators $a$ and $a^\dag$ are unbounded by (\ref{aa}), so,
 the right hand side of (\ref{JC})
is not Lipschitz continuous in $\rho(t)$.
  Accordingly, the meaning of the equation (\ref{JC}) must be adjusted (see Definition \ref{gs}).

In the case of time-independent pumping $A^e(t)=A^e$, the Hamilton 
$H(t)=H$  is a  selfadjoint operator, so
for $\ga=0$,
 solutions are given by
  \be\la{rot}
  \rho(t)=e^{-i Ht}\rho(0)e^{iHt},\qquad t\ge 0.
        \ee
        In this case, the trace $\tr\!\rho(t)$ is conserved, and $\rho(t)\ge 0$ if $\rho(0)\ge 0$.
    However, for time-dependent pumping or $\ga>0$, the formula for solutions is missing.
    \smallskip
    
 We    assume that 
the pumping $A^e(t)$ and 
 the dissipation operator $D(t)$
 are polynomials   in $a$ and $a^\dag$. Moreover, the structure of $D(t)$ corresponds to 
 the theory of
completely positive and trace preserving generators (CPTP) developed by 
Lindblad \ci{L1976} and Gorini, Kossakowski and Sudarshan \ci{GKS1976}.
\smallskip\\
{\bf Examples.}
$A^e=a^\dag+a$, $A^e=a^\dag a$, and  $D(t)=D_1$  satisfy all our assumptions.

 Our main results are as follows:
   there exist 
  continuous linear operators  $U(t):\HS\to \HS$ with $t\ge 0$
 such that for each nonnegative $\rho^0\in\HS$, the trajectory $\rho(t)=U(t)\rho^0$ is
a nonnegative generalised solution
 to  (\ref{JC}) with initial condition
 \be\la{inicw}
 \rho(0)=\rho^0,
 \ee
  and 
    the  a priori bounds hold,
 \be\la{HS3}
\Vert\rho(t)\Vert_\HS\le \Vert\rho^0\Vert_{\HS},\qquad t\ge 0.
\ee
     Our strategy  relies on
  finite-dimensio\-nal approximations  of the annihilation and creation ope\-rators
 and the uniform bounds  (\ref{HS3}).
 The bounds are due to the nonpositivity of the dissipation operators $D(t)$
 under our assumpions.
The key example of such dissipation operators is $D_1$
for which the nonpositivity has been established in \ci{KK2025-S}.

Note that the nonnegativity of $\rho(t)$ cannot be obtained by a straightforward
application of  
 the theory \ci{GKS1976,L1976} of
completely positive and trace preserving generators (CPTP) since the generator $\aA(t)$
is an unbounded operator. However, the theory allows us to
 establish the nonnegativity for special finite-dimensional approximations
 which 
keep the CPTP structure of the generator. Then the nonnegativity of $\rho(t)$ follows by the
 limit transition. In contrast, the trace conservation does not follow in the limit, 
 though it is conserved for all the approximations.

\smallskip

Let us comment on previous related results.
In the case of bounded generators $\aA(t)$,
global solutions for  equations of type (\ref{JC}) obviously exist.
In this case,
Lindblad \cite{L1976} and  Gorini, Kossakowski, and Sudarshan \cite{GKS1976}
found necessary and sufficient conditions on $\aA(t)$ providing the positivity and trace preservation.

For unbounded generators, the existence of global solutions
to Quantum Dynamical Systems (QDS) is not well-developed, \ci[p.110]{AL2007}.
In \ci{D1977}, E.B. Davies considered autonomous quantum-mechanical Fokker--Planck 
equations (QFP).  
The existence of the corresponding positive 
contraction semigroup is established 
  in the Banach space of self-adjoint trace-class operators.
  The uniqueness and trace preservation have not been proved.
  The useful sufficient conditions, providing the trace preservation 
  for QFP equations,
  were 
  found in \ci{CF1998}. The detailed characterisation of a class of covariant 
  QDS with unbounded generators is presented in \ci{H1996}. 
  The survey
  can be found in \ci{D1976,F1999}.

In \ci{KK2025-S}, we have constructed the contraction semigroup for the equation
(\ref{JC}) with 
 pumping  and dissipation operator which do not depend on time, so $\aA(t)=\aA$. In this case, the main ingredient in the proof 
  is the nonpositivity $\aA\le 0$ which implies the existence of 
  the semigroup $e^{\aA t}$ by
   the  Lumer--Phillips theorem.
The nonpositivity
  follows from the same property
  of the dissipation operator $D$ and
  implies that the semigroup 
 is contracting, i.e., the bounds (\ref{HS3}) hold. As the key
example, the nonpositivity  is proved for the operator (\ref{tiGa}). 
   
 Note that
the equation (\ref{JC}) is non-autonomous, so, the  theory of semigroup 
 is not applicable.
 Our approach 
 allows us to substitute the theory  by a systematic application of
 the contraction (\ref{HS3}).

Up to our knowledge,
the well-posedness
for nonautonomous damped driven Jaynes--Cummings equations
has not been
established previously.
\smallskip

{\bf Acknowledgements.} The authors thank S. Kuksin, M.I. Petelin, A. Shnirelman and H. Spohn
 for longterm fruitful discussions, and the  Institute of Mathematics of BOKU University
for the support and hospitality. The research was funded in whole  by
Austrian Science Fund (FWF) 10.55776/PAT3476224.

    \section{Notations and main results}
  
    Every density operator $\rho\in\HS$
    is defined uniquely by its matrix entries
      \be\la{roent}
      \rho_{n,s;n',s'}=\langle n,s|\rho|n',s'\rangle, \qquad n,n'=0,1,\dots,\quad{s,s'=s_\pm}.
      \ee
The Hilbert--Schmidt norm, corresponding to  the inner product  (\ref{HS}), can be written as
\be\la{HSen}
\Vert\rho \Vert_\HS^2=\tr[\rho^2]=\sum_{n,n'=0}^\infty
\sum_{s,s'=s_\pm}|\rho_{n,s;n',s'}|^2<\infty.
\ee
Note that 
by (\ref{aa}), 
the space
$X_\infty$ is invariant with respect to $a$ and $a^\dag$. Hence,
the products 
of
$\rho\in\cD$ with any polynomials of $a$ and $a^\dag$ 
are well-defined as operators
in $X_\infty$.

 We assume that the pumping and the dissipation operator 
satisfy the following conditions. Denote by $M_2=\C^2\otimes\C^2$ the space of 
$2\times 2$-matrices.
\smallskip\\
{\bf H1.} 
 The pumping $A^e(t)$ is a
 symmetric operator  in  $X_\infty$, and it is a
polynomial in the creation and annihilation operators:
 \be\la{poldiA}
 A^e(t)=\sum A^e_{l,l'}(t)a^{\dag l}a^{l'},\qquad A^e_{l,l'}\!\!\in\! C(0,\infty;M_2),
 \ee
{\bf  H2.}
 The dissipation operator $D(t)\rho$
 admits a  ``double" structure
 of the
 CPTP theory 
 of completely positive trace preserving operators \ci{GKS1976,L1976,S1980}:
   \be\la{poldiD}
D(t)\rho=[V(t)\rho, V^\dag(t)]+[V(t),\rho V^\dag(t)],
 \ee
where $V(t)$ is also a
polynomial in the creation and annihilation operators:
\be\la{poldiV}
 V(t)=\sum V_{l,l'}(t)a^{\dag l}a^{l'},\qquad V_{l,l'}\!\!\in\! C(0,\infty;M_2),
\ee

\br\rm
i) {\bf H1--H2} imply that
$\aA(t)\rho\in\cD$ for 
$\rho\in\cD$ and $t\ge 0$.
\smallskip\\
ii) 
The operator (\ref{tiGa}) admits the structure (\ref{poldiV}).  

\er

\br\rm
In \ci{JA1993}, the pumping 
\be\la{Abnp}
A^e(t)=\si_+e^{-i\om_p t} +\si_-e^{i\om_p t},\qquad \si_+=\left(\!\!\!\ba{cc}0&1\\0&0\ea\!\!\!\right)\!,\quad \si_-=\left(\!\!\!\ba{cc}0&0\\1&0\ea\!\!\!\right)
\ee
 has been applied
 to the study of the collapse and revival
phenomenon (see \ci[(2.1)]{JA1993}).

\er


 
To formulate our results, we need 
give a meaning  to
the equation
(\ref{JC}). The issue is that the operators $\aA(t)$ are not well-defined on $\HS$.
The equation admits the treatment via the matrix 
entries (\ref{roent}) as the system
 \be\la{JCw} 
 \dot\rho_{n,s;n',s'}(t)=[\aA(t)\rho(t)]_{n,s;n',s'},
 \qquad t\ge 0,\qquad n,n'\ge 0,\quad s,s'=s_\pm,
 \ee
since the right hand side is well-defined for all $\rho(t)\in\HS$.
Indeed,
by
{\bf H1} and (\ref{aa}),
 the operators $\aA(t)$ are well-defined  on the domain $\cD$, and
 \be\la{JCwen} 
 [\aA(t)\rho]_{n,s;n',s'}=
\sum_{\scriptsize
\ba{c}|k-n|+|k'-n'|\le N\\ r,r'=s_\pm\ea} \aA_{n,s;n',s'}^{k,r;k',r'}(t)\rho_{k,r;k',r'},\,\quad n,n'\ge 0,  s,s'=s_\pm,\quad\rho\in\cD,
 \ee
 where $N=\max(2,\deg A^e, \deg V)$, and
\be\la{JCw2C}
\aA_{n,s;n',s'}^{k,r;k',r'}(\cdot)\in C(0,\infty),\qquad \forall n,s,n',s',k,r,k',r'.
\ee
Finally, since
the  summation in (\ref{JCwen}) is finite,
the matrix entries  admit a unique extension by continuity from $\rho\in\cD$
to all $\rho\in\HS$. The structure (\ref{JCwen}) means that the matrix of the generator in the basis $|n,s\rangle$ is almost
diagonal.

Let us denote by
 $\HS_w$ the Hilbert space $\HS$ endowed with 
the weak topology,
and the space of trajectories 
\be\la{TT}
\cT=L^\infty(0,\infty;\HS)\cap
C(0,\infty;\HS_w).
\ee

\bd\la{gs}
A trajectory $\rho(t)\in 
\cT$ is a generalised solution to (\ref{JC}) 
if it 
satisfies the system (\ref{JCw}) in the sense of distributions, i.e.,
  \be\la{JCw2d} 
\rho_{n,s;n',s'}(t)\!-\!\rho_{n,s;n',s'}(0)
=
\int_0^t
[\aA(\tau)\rho(\tau)]_{n,s;n',s'}d\tau.
%
 \ee

\ed

Our main  result is the following theorem.
Denote by $\HS^+$ the set of nonnegative $\rho\in\HS$.
\bt\la{tm}
Let  conditions {\bf H1--H2}  hold.
  Then 
  there exist 
  continuous linear operators  \linebreak $U(t):\HS\to \HS$ with $t\ge 0$
 such that 
 \smallskip\\
i) For each $\rho^0\in\HS$, the trajectory $\rho(t)=U(t)\rho^0\in\cT$ is
the generalised solution
 to  (\ref{JC}) with initial condition (\ref{inicw});
 \smallskip\\
ii)  $\rho(t)\in\HS^+$ for $t\ge 0$ if  $\rho^0\in\HS^+$,
  and 
 the bounds (\ref{HS3}) hold.

\et
The proof relies on finite-dimensional
Faedo--Galerkin approximations defined via 
finite-dimen\-sio\-nal approximations
 of the annihilation and creation operators,
and on uniform a priori bounds in the Hilbert--Schmidt norm.

\br\rm
i)
Theorem \ref{tm} also holds for dissipation operators $D(t)=\sum_j D_j(t)$ 
where each $D_j(t)$ admits the structure (\ref{poldiD})--(\ref{poldiV}) with the corresponding
polyniomial $V_j(t)$. The proof is unchanged. 
\smallskip\\
ii) The trace conservation holds for every finite-dimensional approximation,
however, we cannot conclude the conservation in the limit.

\er

\setcounter{equation}{0}
\section{Nonpositivity of the dissipation operator}

In this section we prove
the nonpositivity of the dissipation operator $D(t)$ 
satisfying {\bf H2}.

\bl\la{lnp}The operators $D(t)$ with the structure (\ref{poldiD}), (\ref{poldiV}) are nonpositive in $\cD$:
   \be\la{np}
\langle\rho,D(t)\rho\rangle_\HS\le 0,\qquad\rho\in\cD,\quad t\ge 0.
\ee 
\el
\bpr We will omit the dependence on $t$.
For $\rho\in\cD$,
\beqn\la{dtr2}
\langle\rho,D\rho\rangle_\HS&=&\tr\big(\rho D\rho\big)
=\tr\Big(\rho\big(
V\rho V^\dag-\fr12 V^\dag V\rho-\fr12 \rho V^\dag V
\big)
\Big)
\nonumber\\
\nonumber\\
&=&
\tr\big(\rho V\rho V^\dag-\rho V^\dag V\rho
\big)
=\tr\big(\rho V\rho V^\dag- V^\dag V\rho^2
\big).
\eeqn
Now we use the fact that $\rho$ is a finite rank Hermitian operator (\ref{fr1}).
Then (\ref{poldiD}) and (\ref{poldiV})  imply that
 the operators $\rho V\rho V^\dag$ and $V^\dag V\rho^2$
have only  finite 
number of nonzero entries (\ref{roent}), so their traces 
are well-defined. 
Moreover,
$\rho$ admits a finite spectral resolution  in the orthonormal basis of its eigenvectors  $e_i\in X_\infty$:
\be\la{eigro}
\rho=\sum_{i=1}^\nu\rho_i e_i\otimes e_i.
\ee
In this basis, the entries  $\rho_{ij}=\rho_i\de_{ij}$, and 
the entries 
$V_{jk}=\langle e_j,V e_k\rangle$
and $V^\dag_{kl}=\langle e_k,V^\dag e_l\rangle$
of the operators $V$ and $V^\dag$ are well-defined.
Hence, (\ref{dtr2}) implies,
with summation in repeated indices, 
\beqn\la{dtr3}
\langle\rho,D\rho\rangle_\HS&=&
\rho_i\de_{ij}V_{jk}\rho_k\de_{kl}
V^\dag_{li}
-V^\dag_{kl}V_{lj}\rho_j^2\de_{jk}
%
=
\rho_i V_{ik}\rho_k
V^\dag_{ki}
-V^\dag_{kl}V_{lk}\rho_k^2
\nonumber\\[1ex]
&=&
\rho_i V_{ik}\rho_k
V^\dag_{ki}
-V^\dag_{ki}V_{ik}\rho_k^2
=
V_{ik}V^\dag_{ki}(\rho_i\rho_k
-\rho_k^2
)
\nonumber\\
&=&
\fr12\Big(V_{ik}V^\dag_{ki}(\rho_i\rho_k
-\rho_k^2)+V_{ki}V^\dag_{ik}(\rho_k\rho_i
-\rho_i^2) 
)
=
-\fr12
|V_{ik}|^2(\rho_i-\rho_k)^2
\le 0
\eeqn
since $V^\dag_{ik}=\ov V_{ki}$. Hence,  (\ref{np})  is proved.
\epr

\br\rm
The proof of the nonpositivity 
   essentially depends on the symmetry of  $\rho$.
 \er

\setcounter{equation}{0}
\section{Finite-dimensional approximations
and uniform bounds}
Replace the Hilbert space  $X$ by
its  finite-dimensional subspace 
\be\la{Xnu}
X_\nu=(|n\rangle\otimes s_\pm:n\le \nu),\qquad \nu=1,2,\dots.
\ee
The subspace is invariant under action of the 
annihilation operator $a$, and we define by $a_\nu$ its restriction to $X_\nu$:
\be\la{anu}
a_\nu[ |n\rangle\otimes s]=
[a |n\rangle]\otimes s,\qquad 0\le n\le \nu.
\ee
On the other hand, $X_\nu$ is not invariant under action of the 
creation operator $a^\dag$, so we define $a^\dag_\nu$ as the adjoint to $a_\nu$:
\be\la{anu2}
a_\nu ^\dag[|n\rangle\otimes s]=
\left\{\ba{ll}
[a^\dag |n\rangle]\otimes s,& 0\le n< \nu\\
0,& n= \nu
\ea\right|.
\ee
\br\rm
The definitions (\ref{anu}), (\ref{anu2}) are equivalent to known matrix representations
of $a_{n,n'}$ and $a^\dag_{n,n'}$ restricted to $n,n'\le\nu$.
\er

Denote by $\cD_\nu$ the space of Hermitian operators (\ref{fr1}) with $n,n'\le\nu$,
and define the finite-dimensional approximations of the JC dynamics
(\ref{JC}) by
\be\la{JCnu}
\dot \rho_\nu(t)=\aA_\nu(t)\rho_\nu(t):=-i[H_\nu(t),\rho_\nu(t)]+\ga D_\nu(t)\rho_\nu(t),\qquad t\ge 0,
\ee
where $\rho_\nu(t)$ is a linear operator in $X_\nu$,
and $H_\nu(t)$ and $D_\nu(t)$ are defined as $H(t)$ and $D(t)$ with $a$ and $a^\dag$
replaced by $a_\nu$ and $a^\dag_\nu$ respectively.
  By the conditions {\bf H1--H2} and (\ref{anu})--(\ref{anu2}), 
   $H_\nu(t)$ and $D_\nu(t)\rho_\nu$
   are
    Hermitian operators in $X_\nu$. Hence,
  $\aA_\nu(t)\rho_\nu$
  is also Hermitian, so  the linear equation (\ref{JCnu})
admits a unique global Hermitian  solution $\rho_\nu(t)\in C(0,\infty;\cD_\nu)$ for any Hermitian initial state $\rho_\nu(0)\in \cD_\nu$.
Since $\rho^0\in \HS$,
\be\la{rokl2}
\rho^0=\sum_{n,n'}\sum_{s,s'}\rho_{n,s;n',s'}^0|n,s\rangle\langle n',s'|,\qquad \Vert\rho^0\Vert_\HS=
\sum_{n,n'} \sum_{s,s'} |\rho_{n,s;n',s'}^0|^2<\infty.
\ee
Define the approximate initial state by
\be\la{iniro}
\rho_\nu^0=\sum_{n,n'\le\nu} \sum_{s,s'}\rho_{n,s;n',s'}^0|n,s\rangle\langle n',s'|.
\ee
It is the Hermitian operator since $\rho^0$ is, and
\be\la{iniro2}
\Vert\rho_\nu^0\Vert_\HS\le \Vert\rho^0\Vert_\HS.
\ee
Obviously, $\rho_\nu^0\ge 0$ if $\rho^0\ge 0$.
Denote by $\rho_\nu(t) $ the solution to (\ref{JCnu}) with the initial state $\rho^0_\nu$.

\bl\la{lm}
i) The a priori bounds hold
\be\la{HS3nu}
\Vert\rho_\nu(t)\Vert_\HS\le \Vert\rho^0_\nu\Vert_\HS,\qquad t\ge t_0,
\quad \nu=0,1,\dots.
\ee
ii) $\rho_\nu(t)\ge 0$ for $t\ge 0$ if $\rho^0_\nu\ge 0$.

\el
\bpr
{\it ad i)}
It suffices to prove 
 that
the generator $\aA_\nu(t)$  
is nonpositive in $X_\nu$:
\be\la{npnu}
 \langle \rho,\aA_\nu(t)\rho\rangle_\HS\le 0, \qquad \rho\in X_\nu.
\ee
The generator   can be written as
\be\la{JC2}
\aA_\nu(t) \rho_\nu=K_\nu(t)\rho_\nu+\ga D_\nu(t)\rho_\nu,
\ee
where
\be\la{KDnu}
\!\!\!\!\!\!\!\!\!\!\!\!
K_\nu(t)\rho_\nu=-i[H_\nu(t),\rho_\nu],\qquad  D_\nu(t)\rho_\nu=[V_\nu\rho_\nu, V_\nu^\dag]+[V_\nu,\rho_\nu V_\nu^\dag].
%
\ee
Here $V_\nu=V_\nu(t)$ is obtained from the polynomial $V(t)$ by replacement 
of $a$ and $a^\dag$ by $a_\nu$ and $a_\nu^\dag$, respectively, while $V_\nu^\dag=(V_\nu(t))^\dag$ by (\ref{anu2}).
 For $\rho_\nu\in \cD_\nu$, we have 
\beqn\la{fr2}
\langle\rho_\nu, K_\nu(t)\rho_\nu\rangle_\HS
=-i\tr(\rho_\nu^\dag[H_\nu(t),\rho_\nu])=-i\tr(\rho_\nu[H_\nu(t)\rho_\nu-\rho_\nu H(t)])
=0.
\eeqn
On the other hand, 
\be\la{npnu2}
\langle\rho_\nu,D_\nu(t)\rho_\nu\rangle_\HS\le 0,\qquad\rho_\nu\in\cD_\nu
\ee 
 by Lemma \ref{lnp} applied to the operator $D=D_\nu(t)$.
Therefore, (\ref{npnu}) is proved.
\smallskip\\
{\it ad ii)} The nonnegativity preservation holds since,
by {\bf H1--H2},
the  generator $\aA_\nu(t)$ with every $t\ge 0$ admits the structure of the CPTP theory \ci{GKS1976, L1976}.
Note that the theory addresses the case when $H_\nu(t)$ and $D_\nu(t)$ do not depend on $t$.
However, the continuity of the coefficients in (\ref{poldiA}) and (\ref{poldiV}), 
allows us to approximate the time-dependent coefficients by pice-wise constant ones
on the intervals $\ve n<t<\ve(n+1)$, 
for which the nonnegativity  preservation holds. Finally, the
 Gronwall estimate implies that the
 discrepancy converges to zero
as $\ve\to 0$.
\epr
  \br\rm
The zero quadratic form (\ref{fr2})
means that the vector field $K_\nu(t)\rho_\nu$ 
on the right hand side 
of (\ref{JC2}) 
is orthogonal to $\rho_\nu$
in the space $\HS$. Thus, this term 
in (\ref{JC2}) 
corresponds to rotations in $\HS$.
By (\ref{npnu2}),
the second vector field $\ga D_\nu(t)\rho_\nu$ in (\ref{JC2}) 
is the generator of 
contractions
of $\HS$ which correspond to the quantum spontaneous emission. 

\er

\setcounter{equation}{0}
\section{Passage to limit}
The following lemma implies Theorem \ref{tm}.
\bl\la{llim}
There exists a subsequence $\nu'\to\infty$, such that 
\smallskip\\
i) for all
$\rho^0\in \HS$,
\be\la{rocon}
\rho_{\nu'}(\cdot)\to \rho(\cdot),
\ee
where 
the convergence holds 
in $C(0,\infty;\HS_w)$
with limiting functions
 $$
 \rho(\cdot)\in L^\infty(0,\infty;\HS)\cap C(0,\infty;\HS_w);
 $$
 ii)
 the limiting functions
 are generalised solutions to (\ref{JC}) with the initial condition (\ref{inicw});
 \smallskip\\
 iii) the bounds  (\ref{HS3}) hold;
 \smallskip\\
 iv) the maps $U(t):\rho(0)\mapsto \rho(t)$ are linear  in $\HS$;
  \smallskip\\
 v) $\rho(t)\ge 0$ for $t\ge 0$ if $\rho(0)\ge 0$.
\el

\bpr

Using (\ref{JCwen}), 
the system (\ref{JCnu}) can be written as
\be\la{JCnu2} 
 \dot{\rho}_{\nu; n,s;n',s'}(t)=
\sum_{\scriptsize \ba{c}|n-k|+|n'-k'|\le  N\\ k,k'\le\nu,\quad r,r'=s_\pm\ea} \ti\aA_{n,s;n',s'}^{k,r;k',r'}(t)\rho_{\nu; k,r;k',r'}(t),
\qquad t\ge 0,
\quad n,n'\le \nu,\,\,s,s'=s_\pm,
 \ee
 where by (\ref{anu2}),
 \be\la{AA}
 {\rm either} \quad \ti\aA_{n,s;n',s'}^{k,r;k',r'}(t)=\aA_{n,s;n',s'}^{k,r;k',r'}(t)\qquad{\rm or}
  \quad \ti\aA_{n,s;n',s'}^{k,r;k',r'}(t)=0.
 \ee
 By  (\ref{HS3nu}) and (\ref{iniro2}),
\be\la{HS3ma}
\sum_{n,s;n',s'} |\rho_{\nu;n,s;n',s'}(t)|^2\le \Vert \rho^0\Vert_\HS^2,
\qquad t\ge 0.
\ee
Hence, 
\be\la{HS3ma20}
\sup_{\nu\ge \max(n, n'),\,\,t\ge 0}|\rho_{\nu;n,s;n',s'}(t)|\le \Vert\rho^0\Vert_\HS,\qquad \forall\, n,s,n',s'.
\ee
%
Therefore, (\ref{JCnu2}) and (\ref{JCw2C}) imply that for any $T>0$,
\be\la{HS3ma2}
\sup_{ \scriptsize\ba{c}\nu\ge \max(n, n')\\t\in [0,T]\ea}|\dot\rho_{\nu;n,s;n',s'}(t)|\le C_{n,n'}(\Vert\rho^0\Vert_\HS)
<\infty,\quad \forall\, n,s,n',s'.
\ee
Now let us prove Lemma \ref{llim} for each  fixed $\rho^0\in\HS$. 
 First, by the Arzela--Ascoli theorem,
for any fixed  $n,s,n',s'$,
the sequence
$\{\rho_{\nu;n,s;n',s'}(t):\nu\ge 0\}$, is compact in $C(0,\infty)$.
Hence,
 for a subsequence $\nu'\to\infty$ and all $n,s,n',s'$,
\be\la{conap}
\rho_{\nu';n,s;n',s'}(t)\toC \rho_{n,s;n',s'}(t),\qquad t\ge 0,\quad \nu'\to \infty.
\ee
By the Fatou lemma,
the convergence and the bounds (\ref{HS3ma}), (\ref{HS3ma2}) imply that the limiting entries $\rho_{n,s;n',s'}(t)$ represent a trajectory
$\rho(\cdot)\!\in\! \cT=
L^\infty(0,\infty;\HS)\cap
C(0,\infty;\HS_w)$. 
The
bounds (\ref{HS3ma}) imply 
(\ref{HS3}):
\be\la{line}
 \Vert U(t)\rho(0)\Vert_\HS\le\Vert\rho^0\Vert_\HS,
  \ee
  where $U(t):\rho^0\mapsto\rho(t)$.
Second, 
(\ref{AA}) implies that  the limiting entries $\rho_{n,s;n',s'}(t)$ satisfy the system  (\ref{JCw}) 
in the sense of distributions
since 
the number of summands in (\ref{JCwen})  is finite.
At last,  let us verify the initial condition (\ref{inicw}). It suffices to check that
\be\la{inicn}
\rho_{n,s;n',s'}(t)\to \rho_{n,s;n',s'}^0,\qquad t\to +0,\quad\forall n,s,n',s'.
\ee
This convergence holds 
since
 $\rho_{\nu;n,s;n',s'}^0\to \rho_{n,s;n',s'}^0$
by (\ref{iniro2}), and bounds (\ref{HS3ma2}) are uniform in
$\nu\ge \max(n,n')$. Finally, $\rho(t)\ge 0$ if $\rho(0)\ge 0$ since
$\rho_\nu(t)\ge 0$ by Lemma \ref{lm} ii).

It remains to show that the 
subsequence $\nu'$ in (\ref{conap})
can be chosen the same for all 
initial states
 $\rho^0\in\HS$, and the maps $U(t):\HS\to\HS$ 
 for all $t\ge 0$ 
 are linear and preserving the nonnegativity.
 First, 
 the subsequence can be chosen the same for
 all
 $\rho^0\in \cD_\infty:=\cup_{\nu\ge 0}\cD_\nu$
 with rational matrix entries.
   Then  the limit maps
 $U(t)\rho^0\mapsto \rho(t)$  with all $t\ge 0$
 are linear maps $\cD_\infty\to \HS$ over the field 
 of rational numbers, and preserving the nonnegativity.
 
 Second, the linearity and bounds
 (\ref{line}) imply, 
 that the maps $U(t):\cD_\infty\to \HS$ with all $t\ge 0$ are globally Lipschitz 
 in the metric of $\HS$
  with the Lipschitz constant one.  
 Hence,  by the density arguments,
  the convergence (\ref{conap}) holds
 for all
 $\rho^0\in\HS$. Therefore, all the maps $U(t)$ are linear over $\C$ and preserving the nonnegativity.
\epr

\section{Conflict of interest}
We have no conflict of interest.

 \section{Data availability statement}
The manuscript has no associated data.

\end{document}